\documentclass[english]{scrartcl}
\usepackage[T1]{fontenc}
\usepackage[utf8]{inputenc}
\usepackage{babel}
\usepackage{mathtools}
\usepackage{url}
\usepackage{amsmath}
\usepackage{amsthm}
\usepackage{amssymb}
\usepackage{graphicx}
\usepackage{geometry}
\usepackage[pdfusetitle,
 bookmarks=true,bookmarksnumbered=false,bookmarksopen=false,
 breaklinks=false,pdfborder={0 0 1},backref=false,colorlinks=false]
 {hyperref}

\makeatletter
\usepackage{tikz}

\usepackage{scalerel}
\usepackage{stackengine,wasysym}

\makeatother

\providecommand\theoremname{Theorem}
\theoremstyle{plain}
\newtheorem{thm}{\protect\theoremname}
\begin{document}
\title{Notes for a seminar on vortex layer cascades}
\author{Andrés Laín-Sanclemente}
\date{October-November 2026}
\maketitle
\begin{abstract}
After the seminars I recently gave at New York University, the University
of Chicago, Princeton University and the University of Bonn and after
having been asked by some people to publish these notes, I decided
to upload them to arXiv. In recent years, the method of vortex layer
cascades has successfully been applied to construct finite-time singularities
in multiple fluid equations (3D Euler, Boussinesq, IPM, g-SQG...).
In this talk, I will motivate and describe the method and sketch how
it can be applied to the 2D Boussinesq and 2D non-homogeneous Euler
equations. In both systems, the growth of the vorticity is driven
by an accumulated hysteresis effect triggered by the cyclic deformation
of the density layers, whose dynamics is governed by an infinite chain
of degenerate pendula. Our work also settles the previously open question
about global existence of classical solutions of the (forced) 2D non-homogeneous
Euler equations.
\end{abstract}
The main objective of today's talk is to \textbf{explain} what a \textbf{vortex
layer cascade} is and why this method is \textbf{useful} to look for
\textbf{singularities} in incompressible fluids. The results I will
state are \textbf{joint work with Diego Córdoba and Luis Martínez-Zoroa}.
For the sake of concreteness, we will work with a single example:
the 2D Boussinesq system.

\section{Boussinesq system}

In today's talk we will mainly consider the Boussinesq system
\[
{\displaystyle \begin{aligned}\frac{\partial\rho}{\partial t}+u\cdot\nabla\rho & =f_{\rho},\\
\frac{\partial u}{\partial t}+u\cdot\nabla u & =-\nabla P-\rho\hat{e}_{1}+f_{u},\\
\nabla\cdot u & =0,
\end{aligned}
}
\]
where $\rho,P:\mathbb{R}_{t}\times\mathbb{R}^{2}_{x}\to\mathbb{R}$,
$u:\mathbb{R}_{t}\times\mathbb{R}^{2}_{x}\to\mathbb{R}^{2}$ are the
unknowns, $f_{\rho}:\mathbb{R}_{t}\times\mathbb{R}^{2}_{x}\to\mathbb{R}$,
$f_{u}:\mathbb{R}_{t}\times\mathbb{R}^{2}_{x}\to\mathbb{R}^{2}$ are
the source terms and $\rho$ denotes ``density fluctuations'' (and
can therefore take \textbf{either sign}). The local existence theory
for initial data $\rho_{0},u_{0}\in C^{1,\alpha}$, $\alpha\in\left(0,1\right)$,
$\omega_{0},\nabla\rho_{0}\in L^{q}$, $1<q<2$, was established by
Chae, Kim and Nam \cite{ChaeKimNam1999}. Taking the curl in the velocity
equation, one is able to remove the pressure, obtaining the vorticity
equation
\[
{\displaystyle \frac{\partial\omega}{\partial t}+u\cdot\nabla\omega=\frac{\partial\rho}{\partial x_{2}}}+f_{\omega},\quad\omega=\nabla^{\perp}\cdot u,\quad f_{\omega}=\nabla^{\perp}\cdot f_{u},\quad\nabla^{\perp}=\left(-\frac{\partial}{\partial x_{2}},\frac{\partial}{\partial x_{1}}\right).
\]

\section{Strategy: cascade of vorticity layers}

\subsection{Philosophical motivation}

\begin{figure}[tbph]
\begin{centering}
\includegraphics[width=0.8\linewidth]{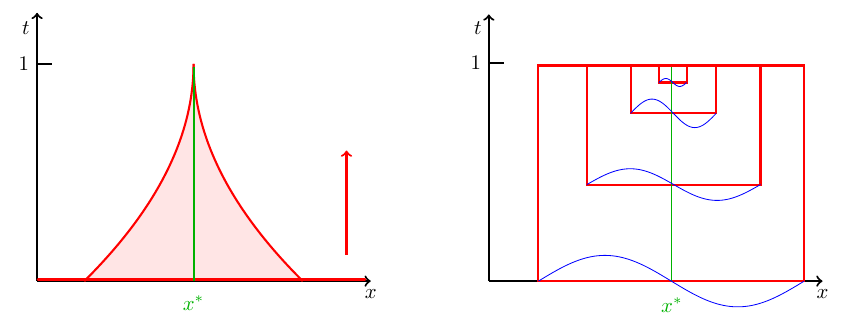}
\par\end{centering}
\caption{\textbf{Left:} self-similar approach. True self-similar solutions
provide \textquotedblleft continuous\textquotedblright{} control over
the whole space-time, while asymptotically self-similar solutions
provide \textquotedblleft continuous\textquotedblright{} control over
a \textquotedblleft cone\textquotedblright{} in space-time, with its
vertex at the blow-up point and time. \textbf{Right:} vortex layer
cascade. The behavior of the solution is only known through a sequence
of triples of spatial supports, spatial frequencies and temporal supports.
The horizontal width of each rectangle denotes its spatial support,
its height represents its temporal support and its blue sine wave
stands for its frequency. As one can see, this method provides control
over a countable sequence of frequencies at the blow-up point and
time.}\label{fig:comparison strategies}

\end{figure}

In order to prove blow-up, one needs a counterexample to global existence,
i.e., an actual solution of the equation that leaves the local-existence
regime at finite time. Explicit solutions of fluid equations are generally
quite rare and the ones we know do not serve for this purpose. Instead
of knowing the whole solution, one has to \textbf{settle for} refined
\textbf{control at infinitesimally small spatial and temporal scales}
near the blow-up point and time. (If the control just reaches finitely
small spatial and temporal scales, only growth can be proven). As
far as I know, there are two main ways of achieving this local control
around the blow-up.
\begin{itemize}
\item \textbf{Self-similar} (or asymptotically self-similar) \textbf{solutions}.
The blow-up is described through a ``static'' profile. With this
method, solving the fluid equation \textbf{reduces to an ODE or a
PDE on a smaller-dimensional space}. {[}Elgindi, Jeong, Chen, Hou,
Pasqualotto, Gómez-Serrano, Cao-Labora, Buckmaster, Lai, Wang, Shao,
Wei, Zhang, Zheng, ...{]}. \cite{Elgindi2021,ElgindiGhoulMasmoudi2019,ElgindiJeong2019,ElgindiJeong2020,ChenHou2021,ChenHouI,ChenHouII,ChenHouHuang2022,Chen2026I,Chen2026II,ElgindiPasqualotto2023,WangLaiGomezSerranoBuckmaster2023,WangEtAl2025,WangLegerLaiBuckmaster2025,ShaoWeiZhangZhang2026,Zheng2023}.
\item \textbf{Vortex layer cascade}. Provides local space-time control around
a (countable) sequence of spatial and temporal scales (i.e., discrete
control). With this method, solving the fluid equation \textbf{reduces
to a countable system of ODEs}. {[}Córdoba, Martínez-Zoroa, Zheng,
Lucas-Manchón, Alpöge, Buckmaster, Dai, OpenAI{]}. \cite{CordobaMartinezZoroaZheng2025,CordobaMartinezZoroa2023,CordobaMartinezZoroa2024,CordobaDominguezLucasManchonMartinezZoroa2026,AlpogeBuckmasterCoiculescu2026IPM,AlpogeBuckmaster2026Boussinesq,AlpogeBuckmaster2026EulerForced,AlpogeBuckmasterCoiculescu2026Euler,Dai2026,OpenAI2026Euler,OpenAI2026NS,Zheng2026IPM}.
\end{itemize}
Figure \ref{fig:comparison strategies} contains a schematic comparison
of both strategies.

\subsection{Implementation}

For concreteness, we will just work with the 2D Boussinesq system
and, most of the time, we will just consider its vorticity equation,
but the same ideas apply to other fluid PDEs.
\begin{description}
\item [{Step1}] \textbf{Find a growth mechanism for the critical norms
of the unknowns.} For example, in the Boussinesq setting, the critical
norms are $\left|\left|\rho\right|\right|_{\dot{C}^{1}}$ and $\left|\left|\omega\right|\right|_{L^{\infty}}$.
As the density equation is a transport equation, growth can only be
achieved through deformation. Then, in view of the vorticity equation,
deformation of the density can lead to an increase in amplitude of
the vorticity. \textbf{Deformation of $\rho$ is determined by $\nabla u$}.
Deformation can happen in multiple ways. Conceptually, it helps to
consider the decomposition of $\nabla u$ given in Figure \ref{fig:deformation gradu}.

\begin{figure}[tbph]
\begin{centering}
\includegraphics[width=1\linewidth]{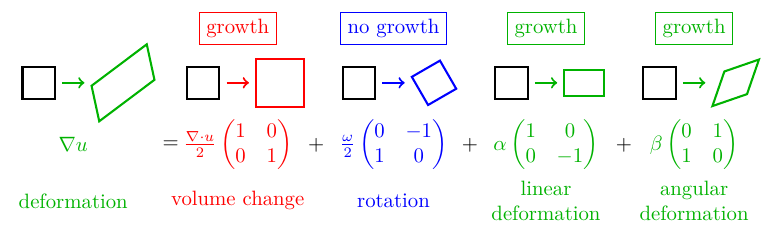}\caption{Schematic decomposition of a general deformation of a 2D fluid volume.
The \textquotedblleft growth\textquotedblright{} and \textquotedblleft no
growth\textquotedblright{} labels refer to the possibility of growth
of $\left|\left|\rho\right|\right|_{\dot{C}^{1}}$. In our case, the
fluid is incompressible, so no volume change is possible. The growth
of our construction for 2D Boussinesq comes from linear deformation.
However, in the construction for IPM and (and all that inherit from
it), the velocity field is a shear flow, which causes a deformation
that is a mixture of rotation and angular deformation. In particular,
the diagonal components of $\nabla u$ are zero.}\label{fig:deformation gradu}
\par\end{centering}
\end{figure}

\item [{Step2}] \textbf{Write the unknowns and source terms as a sum of
layers}. In the Boussinesq setting,
\[
\rho=\sum_{n\in\mathbb{N}}\rho^{\left(n\right)},\quad u=\sum_{n\in\mathbb{N}}u^{\left(n\right)},\quad\omega=\sum_{n\in\mathbb{N}}\omega^{\left(n\right)},\quad f_{\rho}=\sum_{n\in\mathbb{N}}f^{\left(n\right)}_{\rho},\quad f_{\omega}=\sum_{n\in\mathbb{N}}f^{\left(n\right)}_{\omega}.
\]
In what follows, we will use capital letters to denote partial sums
up to layer $n$:
\[
P^{\left(n\right)}=\sum^{n}_{m=1}\rho^{\left(m\right)},\quad U^{\left(n\right)}=\sum^{n}_{m=1}u^{\left(m\right)},\quad\Omega^{\left(n\right)}=\sum^{n}_{m=1}\omega^{\left(m\right)},\quad F^{\left(n\right)}_{\rho}=\sum^{n}_{m=1}f^{\left(m\right)}_{\rho},\quad F^{\left(n\right)}_{\omega}=\sum^{n}_{m=1}f^{\left(m\right)}_{\omega}.
\]
\item [{Step3}] \textbf{Substitute Step2 into the fluid equations and impose
that the partial sums up to layer $n$ be a solution of the equations
$\forall n\in\mathbb{N}$}. By subtracting the equation that the partial
sums up to layer $n-1$ satisfy from the equation that the partial
sums up to layer $n$ satisfy, we obtain the $n$-th layer equation,
i.e., the equation that the unknowns of layer $n$ have to fulfill.
In our case, this amounts to enforcing that $\left(P^{\left(n\right)},\Omega^{\left(n\right)},F^{\left(n\right)}_{\rho},F^{\left(n\right)}_{\omega}\right)$
must solve the 2D Boussinesq equation $\forall n\in\mathbb{N}$, which
provides layer equations for $\left(\rho^{\left(n\right)},\omega^{\left(n\right)},f^{\left(n\right)}_{\rho},f^{\left(n\right)}_{\omega}\right)$
$\forall n\in\mathbb{N}$. The vorticity equation of layer $n$ becomes
\[
\frac{\partial\omega^{\left(n\right)}}{\partial t}+U^{\left(n-1\right)}\cdot\nabla\omega^{\left(n\right)}+u^{\left(n\right)}\cdot\nabla\Omega^{\left(n-1\right)}+u^{\left(n\right)}\cdot\nabla\omega^{\left(n\right)}=\frac{\partial\rho^{\left(n\right)}}{\partial x_{2}}+f^{\left(n\right)}_{\omega}.
\]
In general, a \textbf{countable triangular system of PDEs} is obtained.
Moreover, terms that were linear in the unknowns in the original equation
will involve only layer $n$ in the $n$-th layer equation; whereas
quadratic terms in the original equations will involve not only layer
$n$ but also all layers that precede it. Up to this point, \textbf{no
real simplification} of the original fluid PDE has been achieved.
Quite the contrary, \textbf{each layer equation is at least as hard
to solve as the original fluid equation!}
\item [{Step4}] \textbf{Assume that the pure quadratic term} (the only
one involving two quantities associated to the same layer) \textbf{is
negligible} with respect to the linear terms. (I will comment more
on how to achieve this in Step8). If this is the case, the force may
be used to absorb the errors committed in the approximation. Then,
the $n$-th layer equation becomes \textbf{linear in the unknowns!}
In the case of the Boussinesq system, we obtain
\[
\frac{\partial\omega^{\left(n\right)}}{\partial t}+U^{\left(n-1\right)}\cdot\nabla\omega^{\left(n\right)}+u^{\left(n\right)}\cdot\nabla\Omega^{\left(n-1\right)}\approx\frac{\partial\rho^{\left(n\right)}}{\partial x_{2}}+f^{\left(n\right)}_{\omega}.
\]
\item [{Step5}] (Optional) \textbf{Exploit the different spatial and temporal
scales to neglect more terms of the equation}. For example, in the
Boussinesq case, since the critical norm is $\left|\left|\omega\right|\right|_{L^{\infty}}$,
it should be possible for $\left|\left|u\right|\right|_{L^{\infty}}$
to remain bounded at the blow-up time. In that case, we would have
$\left|\left|u^{\left(n\right)}\right|\right|_{L^{\infty}}\xrightarrow[n\to\infty]{}0$.
Suppose $u^{\left(n\right)}$ has amplitude $A_{n}$ and frequency
$q_{n}$. Then, $\nabla\omega^{\left(n\right)}\sim A_{n}q^{2}_{n}$,
so
\[
\left|\left|u^{\left(n\right)}\cdot\nabla\Omega^{\left(n-1\right)}\right|\right|_{L^{\infty}}\sim A_{n}\sum^{n-1}_{m=1}A_{m}q^{2}_{m}.
\]
If $A_{m}$ decreases superexponentially in $m\in\mathbb{N}$ (for
example, $A_{m}=\mathrm{e}^{-\mathrm{e}^{m}}$) and $A_{m}q^{2}_{m}$
grows superexponentially in $m\in\mathbb{N}$ (for example, $A_{m}q^{2}_{m}=\mathrm{e}^{\mathrm{e}^{m}}$),
the main order of the sum above will be given by its last summand,
i.e.,
\[
\left|\left|u^{\left(n\right)}\cdot\nabla\Omega^{\left(n-1\right)}\right|\right|_{L^{\infty}}\sim A_{n}A_{n-1}q^{2}_{n-1}
\]
and the amplitude $A_{n}$ will be able to compensate the growth of
$A_{n-1}q^{2}_{n-1}$. In conclusion, by choosing the amplitudes and
frequencies in this way, we may also absorb this term with the force.
Thence, we arrive to
\[
\frac{\partial\omega^{\left(n\right)}}{\partial t}+U^{\left(n-1\right)}\cdot\nabla\omega^{\left(n\right)}\approx\frac{\partial\rho^{\left(n\right)}}{\partial x_{2}}+f^{\left(n\right)}_{\omega}.
\]
\item [{Step6}] \textbf{Neglecting the force, solve the linear transport
equation that arises.} \textbf{Optionally}, one may \textbf{exploit
the different spatial and temporal scales to replace the velocity
of previous layers by its first-order Taylor approximation}, making
the transport equation much easier to solve. The idea is the following:
if the support of layer $n$ is going to be much smaller than the
support of layer $n-1$, it makes sense that just a first-order approximation
of $U^{\left(n-1\right)}$ is needed to understand the evolution of
layer $n$. In the case of Boussinesq, by Taylor's Remainder Theorem,
the error we make in this approximation is of order
\[
\left|\left|\mathrm{D}^{2}U^{\left(n-1\right)}\right|\right|_{L^{\infty}}\text{diam}\left(\text{supp}\left(\omega^{\left(n\right)}\right)\right)^{2}\left|\left|\nabla\omega^{\left(n\right)}\right|\right|_{L^{\infty}}\sim\sum^{n-1}_{m=1}A_{m}q^{2}_{m}\frac{1}{q^{2}_{n}}A_{n}q^{2}_{n}\sim A_{n-1}q^{2}_{n-1}A_{n}.
\]
(This computation is not exactly true because $\text{diam}\left(\text{supp}\left(\omega^{\left(n\right)}\right)\right)$
does not actually decrease as fast as $\frac{1}{q_{n}}$, but it is
a good enough approximation to understand the main idea). We obtain
exactly the same as in Step5. Therefore, the error is controlled and
can be absorbed into the force. 

If $\Phi^{\left(n\right)}$ is the true flow map of layer $n$,
\[
\frac{\partial\Phi^{\left(n\right)}}{\partial t}\left(t,x\right)=U^{\left(n-1\right)}\left(t,\Phi^{\left(n\right)}\left(t,x\right)\right),\quad\Phi^{\left(n\right)}\left(0,x\right)=x,
\]
the \textbf{first-order approximate flow-map} $\phi^{\left(n\right)}\left(t,x\right)=\eta^{\left(n\right)}\left(t\right)+\left(x-c^{\left(n\right)}\right)M^{\left(n\right)}\left(t\right)$
is obtained by solving
\[
{\small \underbrace{\frac{\partial\eta^{\left(n\right)}}{\partial t}\left(t\right)=U^{\left(n-1\right)}\left(t,\eta^{\left(n\right)}\left(t\right)\right),\quad\eta^{\left(n\right)}\left(0\right)=c^{\left(n\right)}}_{\text{displacement of the layer}},\quad\underbrace{\frac{\partial M^{\left(n\right)}}{\partial t}\left(t\right)=M^{\left(n\right)}\left(t\right)\cdot\nabla U^{\left(n-1\right)}\left(t,\eta^{\left(n\right)}\left(t\right)\right),\quad M^{\left(n\right)}\left(0\right)=\text{id}}_{\text{deformation of the layer}},}
\]
where $c^{\left(n\right)}$ denotes the initial position of the layer
center. In particular, the equation for the deformation of the layer
can be integrated to
\begin{equation}
M^{\left(n\right)}\left(t\right)=\mathcal{T}\exp\left(\int^{t}_{0}\nabla U^{\left(n-1\right)}\left(s,\eta^{\left(n\right)}\left(s\right)\right)\mathrm{d}s\right),\label{eq:M(t)}
\end{equation}
where $\mathcal{T}\exp$ denotes a time-ordered exponential. This
last equation yields one key insight: If layer $n$ has separation
of spatial variables initially and we want this separation of spatial
variables to prevail in time, $\nabla U^{\left(n-1\right)}\left(t,\eta^{\left(n\right)}\left(t\right)\right)$
must be diagonal. In any other case, some shear will inevitably occur.

In this manner, we have managed to describe the transport and deformation
of each layer through two ``very simple'' ODEs. Moreover, this procedure
allows us to approximately eliminate the transport term from the layer
equation. In the case of Boussinesq, we arrive at
\begin{equation}
\frac{\partial}{\partial t}\left(\omega^{\left(n\right)}\left(t,\phi^{\left(n\right)}\left(t,x\right)\right)\right)\approx\frac{\partial\rho^{\left(n\right)}}{\partial x_{2}}\left(t,\phi^{\left(n\right)}\left(t,x\right)\right),\label{eq:growth mechanism}
\end{equation}
the equation that describes the growth mechanism for the vorticity.
If one were now able to use separation of variables between space
and time and to choose the spatial dependence of $\omega^{\left(n\right)}\left(t,\phi^{\left(n\right)}\left(t,x\right)\right)$
and $\frac{\partial\rho^{\left(n\right)}}{\partial x_{2}}\left(t,\phi^{\left(n\right)}\left(t,x\right)\right)$
to be the same, the $n$-th layer equation reduces to an ODE which
describes the growth of the amplitude of layer $n$ in time.

To sum up, in this last step, we have successfully simplified the
layer PDE to a finite system of ODEs describing the \textbf{displacement,
deformation and growth of each layer}. In total, we have \textbf{replaced
the original fluid PDE by a countable triangular system of ODEs!}
\item [{Step7}] Now that the displacement, deformation and growth mechanism
are known, depending on the construction, there may still be some
leeway to choose the characteristic time scales of the layers. If
the \textbf{time scales can be chosen decreasing sufficiently fast}
in $n\in\mathbb{N}$, a \textbf{cascade in the form $1\to2\to3\to\dots\to n\to\dots$
can be established}. Then, the layers are introduced sequentially
in time and the growth occurs as follows: 1. Layer $n$ is created
through the force. 2. Layer $n-1$ makes layer $n$ grow. 3. Layer
$n+1$ is introduced through the force. 4. Layer $n$ makes layer
$n+1$ grow, etc. Moreover, since the evolution time of layer $n$
is much smaller than the evolution time of layer $n-1$, for the purposes
of layer $n$, layer $n-1$ can be considered \textbf{constant in
time}, as if it had already attained the state it will have at the
blow-up time.
\item [{Step8}] \textbf{Choose the functional form of the unknowns so that
the pure quadratic term} (the only one involving two quantities associated
to the same layer) \textbf{becomes negligible} with respect to the
linear terms (so as to justify Step4). There is\textbf{ no standard
way }to carry out this step, which is highly dependent on the equation
considered. Nonetheless, the following are usually required:
\begin{enumerate}
\item The original fluid equation must have a \textbf{family of stationary
solutions}. Then, these non-zero stationary solutions are candidates
for the functional form, as they encode the required cancellation.
\item The \textbf{family of stationary solutions must be closed under transport
by} $U^{\left(n-1\right)}$, i.e., transporting any member of the
family by $U^{\left(n-1\right)}$ provides another stationary solution
of the same family. This is related to the quantity $M^{\left(n\right)}\left(t\right)$
(see equation (\ref{eq:M(t)})) and usually requires $\nabla U^{\left(n-1\right)}$
to obey a certain type of symmetry. In the Boussinesq case, $M^{\left(n\right)}\left(t\right)$
will be diagonal.
\item The family must be \textbf{compatible with the growth mechanism}.
For instance, in the case of our Boussinesq construction, it must
allow for separation of temporal and spatial variables in (\ref{eq:growth mechanism})
and to choose $\omega^{\left(n\right)}\left(t,\phi^{\left(n\right)}\left(t,x\right)\right)$
and $\frac{\partial\rho^{\left(n\right)}}{\partial x_{2}}\left(t,\phi^{\left(n\right)}\left(t,x\right)\right)$
with the same spatial dependence.
\end{enumerate}
In the case of Boussinesq, if one aims for linear deformation, it
can be shown that trigonometric functions are the only ``reasonable''
building blocks that satisfy these three requirements. However, trigonometric
functions do not have compact support. Fortunately, the product of
a sufficiently ``wide'' cut-off function and trigonometric functions
is an ``almost''-stationary solution, so that the errors are small
enough to be absorbed by the force.
\item [{Step9}] (Usually optional) \textbf{Improve the regularity of the
force}. If the force resulting from the procedure above is regular
enough and has an adequate functional form, its regularity can be
further improved by adding sublayers specially engineered to cancel
out the least-regular terms of the force and by increasing the order
of the Taylor approximation utilized for the approximate transport.
The idea is that, as the force is regular enough, the sublayers needed
to eliminate the least-regular term of the force will have amplitudes
and frequencies much smaller than those of the main layers of the
construction, allowing us to treat the sublayers as perturbations
of the layers. Iterating this procedure, one can make the force $C^{k}$
for any $k\in\mathbb{N}$. If, in addition, the number of correcting
sublayers and the order of the Taylor approximations are chosen to
approach infinity at the blow-up time, the force can be made smooth.
\end{description}
If all simplifying assumptions mentioned in the steps above are met,
the evolution of a layer reduces to studying the properties described
in Figure \ref{fig:summary}.

\begin{figure}[tbph]
\begin{centering}
\includegraphics[width=0.8\linewidth]{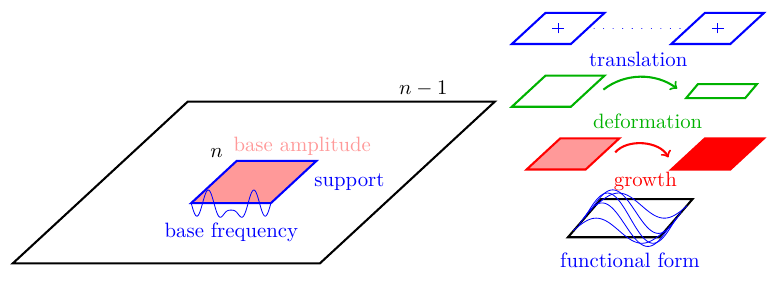}
\par\end{centering}
\caption{\textbf{}Summary of the most important time-independent and time-dependent
properties that define the evolution of a layer. }\label{fig:summary}
\end{figure}

\section{Boussinesq construction}

Unfortunately, I do not have the time to talk about the details of
the Boussinesq construction. Let me just make a few remarks:
\begin{itemize}
\item The velocity field of each layer contains \textbf{hyperbolic points},
which cause the linear deformation of the next density layer.
\item The ODE for $\eta^{\left(n\right)}_{1}\left(t\right)$, the horizontal
position of the center of layer $n$, is deeply related to that of
a pendulum. Indeed, the ODE for $\Xi^{\left(n\right)}\left(t\right)\coloneqq\eta^{\left(n\right)}_{1}\left(t\right)-\eta^{\left(n-1\right)}_{1}\left(t\right)$,
the horizontal position of the center of layer $n$ with respect to
the horizontal position of layer $n-1$, looks roughly like $\dot{\Xi}^{\left(n\right)}=\sin\Xi^{\left(n\right)}$,
the equation of a \textbf{degenerate inverted half-pendulum!} Indeed,
if we differentiate again with respect to time, $\ddot{\Xi}=\left(\cos\Xi\right)\dot{\Xi}=\cos\Xi\sin\Xi=\frac{1}{2}\sin\left(2\Xi\right)$
and, equivalently, $\frac{\partial^{2}}{\partial t^{2}}\left(2\Xi\right)-\sin\left(2\Xi\right)=0$.
``Inverted'' comes from the minus sign instead of the plus sign
of an ordinary pendulum, half comes from the $2$ in $2\Xi^{\left(n\right)}$
and ``degenerate'' from the fact that, in reality, the second order
ODE comes from a first order ODE.
\item The deformation matrix $M^{\left(n\right)}\left(t\right)$ is \textbf{diagonal}.
This is key to avoid the appearance of shear, i.e., to preserve the
separation of the variables $x_{1}$ and $x_{2}$ of the layers.
\item The vorticity layer equation reduces to the following relation between
the amplitude of the vorticity $A^{\left(n\right)}_{\omega}$, the
amplitude of the density $A^{\left(n\right)}_{\rho}$ and a certain
function $b^{\left(n\right)}$ which measures the deformation of layer
$n$:
\[
A^{\left(n\right)}_{\omega}\left(t\right)=\int^{t}_{0}\underbrace{b^{\left(n\right)}\left(s\right)}_{\geq0}\underbrace{A^{\left(n\right)}_{\rho}\left(s\right)}_{\geq0}\mathrm{d}s.
\]
The \textbf{hysteresis} \textbf{effect} on the amplitude of the vorticity
follows from the equation above.
\end{itemize}
Without further ado, we state our Theorem \cite{CordobaLainMartinezZoroa2025}:
\begin{thm}
Let $\alpha\in\left(0,\alpha_{*}\right)$, $\alpha_{*}=\sqrt{\frac{4}{3}}-1$.
There are classical solutions $\left(u,\rho\right)$ on $\left[0,1\right)\times\mathbb{R}^{2}$
of the forced Boussinesq system such that
\begin{enumerate}
\item There is a \textbf{singularity}:
\[
\lim_{T\to1^{-}}\int^{T}_{0}\left|\left|\nabla\rho\left(t,\cdot\right)\right|\right|_{L^{\infty}\left(\mathbb{R}^{2};\mathbb{R}^{2}\right)}\mathrm{d}t=\infty,
\]
\item The \textbf{forces} remain \textbf{``regular''}:
\[
f_{\rho}\in C^{0}_{t}C^{1,\alpha}_{x,c}\left(\left[0,1\right]\times\mathbb{R}^{2}\right),\quad f_{u}\in C^{0}_{t}C^{1,\alpha}_{x}\left(\left[0,1\right]\times\mathbb{R}^{2}\right)\cap C^{0}_{t}L^{2}_{x}\left(\left[0,1\right]\times\mathbb{R}^{2}\right).
\]
\end{enumerate}
\end{thm}

\section{Non-homogeneous Euler equations}

Before introducing the non-homogeneous Euler equations, let us reflect
briefly on the flexibility of the Boussinesq construction. Because
of how the vortex layer cascade works, far from the blow-up point,
only finitely many layers are non-zero. Hence, we may modify the solution
\textbf{far} from the blow-up point and compensate the changes through
the forces without sacrificing regularity. Alternatively, we can \textbf{even
modify the equations themselves} as long as the new equations \textbf{still
look like Boussinesq near the blow-up point}. This is basically the
strategy that we will follow. We will find a solution around which
the Boussinesq and non-homogeneous Euler equations look alike at leading
order.

Let me now present the non-homogeneous Euler equations
\[
{\displaystyle \begin{aligned}\frac{\partial\rho}{\partial t}+u\cdot\nabla\rho & =f_{\rho},\\
\frac{\partial u}{\partial t}+u\cdot\nabla u & =-\frac{\nabla P}{\rho}+f_{u},\\
\nabla\cdot u & =0,
\end{aligned}
}
\]
where $\rho,P:\mathbb{R}_{t}\times\mathbb{R}^{2}_{x}\to\mathbb{R},u:\mathbb{R}_{t}\times\mathbb{R}^{2}_{x}\to\mathbb{R}^{2}$
are the unknowns, $f_{\rho}:\mathbb{R}_{t}\times\mathbb{R}^{2}_{x}\to\mathbb{R}$,
$f_{u}:\mathbb{R}_{t}\times\mathbb{R}^{2}_{x}\to\mathbb{R}^{2}$ are
the source terms and $\rho$ denotes a \textbf{physical} density (and
consequently must be non-negative). The local existence theory for
initial data $\rho_{0},u_{0}\in C^{1,\alpha}$, $\alpha\in\left(0,1\right)$,
$u_{0}\in L^{2}$, $\rho_{0}>0$ was established by Danchin and Fanelli
\cite{DanchinFanelli2011}. The dichotomy between global existence
and singularity formation (in any class with local well-posedness)
was open prior to our work \cite{CordobaLainMartinezZoroa2026}. Taking
the curl in the velocity equation, one obtains the vorticity equation
\[
{\displaystyle \frac{\partial\omega}{\partial t}+u\cdot\nabla\omega=\frac{\nabla P}{\rho}\cdot\frac{\nabla^{\perp}\rho}{\rho}+\nabla^{\perp}\cdot f_{u}}.
\]
Now, if we use the velocity equation to solve for $\frac{\nabla P}{\rho}$,
we obtain
\[
\frac{\partial\omega}{\partial t}+u\cdot\nabla\omega=-\left[\frac{\partial u}{\partial t}+u\cdot\nabla u-f_{u}\right]\cdot\frac{\nabla^{\perp}\rho}{\rho}+\nabla^{\perp}\cdot f_{u}.
\]
The key idea to transfer the Boussinesq blow-up to the non-homogeneous
Euler equations is to work with the \textbf{ansatz}
\[
\rho=\overline{\rho}+\rho_{B},\quad u=u_{B},
\]
where $\left(\rho_{B},u_{B},f_{\rho_{B}},f_{u_{B}}\right)$ is one
of our Boussinesq solutions \cite{CordobaLainMartinezZoroa2025} and
$\overline{\rho}$ will be taken large enough. Substituting into the
non-homogeneous Euler equations, we obtain the following system for
the forces $\left(f_{\rho},f_{u}\right)$:
\[
\begin{aligned}f_{\rho_{B}} & =f_{\rho},\\
\frac{\partial\rho_{B}}{\partial x_{2}}\overbrace{\left(1+\frac{1}{\overline{\rho}+\rho_{B}}\left(f_{u,1}-\left(\frac{\partial u_{B}}{\partial t}+u_{B}\cdot\nabla u_{B}\right)_{1}\right)\right)}^{\left(*\right)}+\\
-\frac{1}{\overline{\rho}+\rho_{B}}\frac{\partial\rho_{B}}{\partial x_{1}}\left(f_{u,2}-\left(\frac{\partial u_{B}}{\partial t}+u_{B}\cdot\nabla u_{B}\right)_{2}\right)+f_{\omega_{B}} & =\nabla^{\perp}\cdot f_{u}.
\end{aligned}
\]
The vorticity equation above poses two main difficulties: 1. $f_{u}$
and $\nabla^{\perp}\cdot f_{u}$ appear in the same equation. 2. $\int^{1}_{0}\left|\left|\frac{\partial\rho_{B}}{\partial x_{2}}\left(t,\cdot\right)\right|\right|_{L^{\infty}\left(\mathbb{R}^{2}\right)}\mathrm{d}t=\infty$,
which, a priori, could prevent $f_{u}$ from lying in the well-posedness
regime $f_{u}\in C^{1,\alpha}$, $\alpha>0$, because of the term
that contains $\frac{\partial\rho_{B}}{\partial x_{2}}$. To solve
these issues, we set up a \textbf{fixed-point argument} to compute
$f_{u}$, while exploiting the \textbf{gauge symmetry} between $f_{u}$
and the pressure to make $\left(*\right)$ vanish $\forall t\in\left[0,1\right]$.

To finish, let me present the statement of the result we obtain \cite{CordobaLainMartinezZoroa2026}.
\begin{thm}
Let $\alpha\in\left(0,\alpha_{*}\right)$, $\alpha_{*}=\sqrt{\frac{4}{3}}-1$.
There are classical solutions $\left(u,\rho\right)$ on $\left[0,1\right)\times\mathbb{R}^{2}$
of the forced non-homogeneous Euler equations such that
\begin{enumerate}
\item There is a \textbf{singularity}:
\[
\lim_{T\to1^{-}}\int^{T}_{0}\left|\left|\omega\left(t,\cdot\right)\right|\right|_{L^{\infty}\left(\mathbb{R}^{2}\right)}\mathrm{d}t=\infty,
\]
\item The forces remain \textbf{``regular''}:
\[
f_{\rho}\in C^{0}_{t}C^{1,\alpha}_{x,c}\left(\left[0,1\right]\times\mathbb{R}^{2}\right),\quad f_{u}\in C^{0}_{t}C^{1,\alpha}_{x}\left(\left[0,1\right]\times\mathbb{R}^{2}\right)\cap C^{0}_{t}L^{2}_{x}\left(\left[0,1\right]\times\mathbb{R}^{2}\right).
\]
\end{enumerate}
\end{thm}

\end{document}